\documentclass[A4paper,11pt]{article}
\usepackage[latin1]{inputenc}
\usepackage[T1]{fontenc}
\usepackage{amsfonts,amsmath,stmaryrd}

\title{Pincement sur le spectre et le volume en courbure de Ricci positive}

\begin{document}
\def\R{\rm R}
\def\Re{\mathbb{R}}
\def\ppi{{\boldmath $\pi$}}
\def\S{\mathbb{S}}
\def\N{\mathbb{N}}
\def\sn{\mathbb{S}^{n-1}}
\def\sh{{\rm sh}}
\def\ch{{\rm ch}}
\def\Car{1\hskip-.9mm{\rm l}}
\def\cupp{\mathop{\cup}\limits}
\def\capp{\mathop{\cap}\limits}
\def\Ric{\mathop{\rm Ric}\nolimits}
\def\K{\mathop{\rm K}\nolimits}
\def\lap{\triangle}
\def\supp{\mathop{\rm sup}\limits}
\def\tend{\mathop{\rightarrow}\limits}
\def\inff{\mathop{\rm inf}\limits}
\def\Dim{\mathop{\rm Dim}\nolimits}
\def\dim{\mathop{\rm dim}\nolimits}
\def\Diam{\mathop{\rm Diam}\nolimits}
\def\Ind{\mathop{\rm Ind}\nolimits}
\def\Vol{\mathop{\rm Vol}\nolimits}
\def\card{\mathop{\rm card}\nolimits}
\def\prodd{\mathop{\prod}\limits}
\def\summ{\mathop{\sum}\limits}
\def\limsupp{\mathop{\overline{\lim}}\limits}
\def\liminff{\mathop{\underline{\lim}}\limits}
\def\rn+1{\mathbb{R}^{n+1}}
\def\sn{\mathbb{S}^n}
\def\T{\mathbb{T}}
\def\Ric{\mathop{\rm Ric}\nolimits}
\def\tr{\mathop{\rm tr}\nolimits}
\def\exp{\mathop{\rm exp}\nolimits}
\def\lap{\triangle}
\def\Cut{\mathop{\rm Cut}\nolimits}
\def\Rad{\mathop{\rm Rad}\nolimits}
\def\Diam{\mathop{\rm Diam}\nolimits}
\def\Vol{\mathop{\rm Vol}\nolimits}
\def\card{\mathop{\rm card}\nolimits}
\def\summ{\mathop{\sum}\limits}
\def\dis{\displaystyle}
\def\neg{\hbox{\ni\char'12}}
\def\barint{\neg\kern-7pt\dis\int}
\def\vrc{variété riemannienne compacte }
\def\dis{\displaystyle}
\font\li=linew10 at 8pt
\def\Neg{\hbox{\li\char'12}}
\def\Barint{\Neg\kern-9.5pt\dis\int}
\font\ni=linew10 at 6pt
\def\neg{\hbox{\ni\char'12}}
\def\barint{\neg\kern-7pt\dis\int}
\def\vrc{variété riemannienne compacte}
\def\proof{\noindent{\sl Preuve.~}}
\def\remark{\noindent{\bf Remarque.~}}

\newtheorem{theo}{Théorème}
\newtheorem{proposition}[theo]{Proposition}  
\newtheorem{lemma}[theo]{Lemme}
\newtheorem{example}[theo]{Exemple}
\newtheorem{corollary}[theo]{Corollaire}

\date{}
\maketitle 

\noindent{\sc Erwann AUBRY\footnote{Travaux en partie financés par la bourse FNRS Suisse n°20-101469}}\\[5mm]
\noindent Institut de mathématiques\\
          Université de Neuchâtel,\\
          erwann.aubry@unine.ch

\abstract{\small We shall show that a complete manifold of dimension $n$ with $\Ric\geq n{-}1$ and its $n$-st eigenvalue close to $n$ is both Gromov-Hausdorff close and diffeomorphic to the sphere. This extends, in an optimal way, a result of P.~Petersen \cite{Pe} (as a by-product, we fill a gap stated in the erratum \cite{Pe2}). We shall also show that a manifold with $\Ric\geq n{-}1$ and volume close to $\frac{\Vol\sn}{\#\pi_1(M)}$ is both Gromov-Hausdorff close and diffeomorphic to the (lens) space $\frac{\sn}{\pi_1(M)}$. This extends results of T.~Colding \cite{Co1} and T.~Yamaguchi \cite{Yam2}.\\

Nous montrons qu'une variété complète de dimension $n$, de courbure $\Ric\geq n{-}1$ et dont la n-ième valeur propre est proche de $n$ est Gromov-Hausdorff proche de $(\sn,can)$ et difféomorphe à $\sn$. Ce résultat étend de manière optimale un résultat de P.~Petersen \cite{Pe} (au passage nous comblons le trou annoncé par l'auteur dans l'erratum \cite{Pe2}). Nous montrons également qu'une variété vérifiant $\Ric\geq n{-}1$ et de volume proche de $\frac{\Vol\sn}{\#\pi_1(M)}$ est difféomorphe à l'espace (lenticulaire) $\frac{\sn}{\pi_1(M)}$ et Gromov-Hausdorff proche de la métrique canonique. Ceci améliore des résultats de T.~Colding \cite{Co1} et T.~Yamaguchi \cite{Yam2}.}

\section{Introduction}

{\bf Dans cet l'article, $M$ désigne par défaut une variété riemannienne complète de dimension $n$ et de courbure de Ricci $\Ric\geq(n{-}1)$.} Dans cet ensemble de variétés, la sphère $\sn$ (munie de sa métrique canonique) est un point extrémal pour de nombreux invariants riemanniens. Des théorèmes de comparaison classiques dûs à S. Myers, R.~Bishop, A.~Lichnerowicz, S.~Cheng et M.~Obata astreignent les variétés de cet ensemble à vérifier $\Diam M\leq\Diam\sn=\pi$ (et donc $M$ est compacte et $\pi_1(M)$ est fini), $\Vol M\leq\Vol\sn$, $\Rad M\leq\Rad\sn=\pi$ (où $\Rad$ désigne le radius, i.e. le rayon de la plus petite boule géodésique qui recouvre toute la variété), $\lambda_1(M)\geq\lambda_1(\sn)=n$ (où $0=\lambda_0(M)<\lambda_1(M)\leq\lambda_2(M)\leq\ldots$ désigne le spectre du laplacien sur la variété, compté avec multiplicité). De plus, si l'égalité est réalisée dans une de ces inégalités, alors la variété $M$ est isométrique à la sphère canonique $\sn$ (auquel cas, $n$ est une valeur propre de $M$ de multiplicité $n+1$).

Au vu de ces faits, il est naturel de chercher les propriétés (topolo\-giques, différentiables ou métriques) de $\sn$ qui sont conservées par les variétés riemanniennes de courbure $\Ric\geq n{-}1$ pour lesquelles l'un des invariants riemanniens considérés plus haut prend une valeur suffisamment proche de sa valeur extrémale.

Dans \cite{Co1} et \cite{Co2}, T.~Colding montre que les trois conditions suivantes sont équivalentes (toujours sous l'hypothèse $\Ric\geq n{-}1$):

(1) $\Vol M$ est proche de $\Vol\sn$,

(2) $\Rad M$ est proche de $\pi$,

(3) $M$ est Gromov-Hausdorff proche de $\sn$.

P.~Petersen a par la suite prouvé que la condition $(2)$ est équivalente à ce que les $n{+}1$ premières valeurs propres (non nulles) du laplacien de $M$ soient proches de $n$ (voir \cite{Pe}). De plus, ces 4 conditions impliquent que $M$ est difféomorphe à $\sn$, d'après le:

\begin{theo}[J.~Cheeger-T.~Colding, \cite{Ch-Co3}]\label{Che-Col}
  Soient $K$ un réel fixé et $(M^n_p,g_p)_{p\in\N}$ une suite de variétés riemanniennes complètes vérifiant $\Ric_{g_p}\geq K$ et qui converge au sens de Gromov-Hausdorff vers une variété riemannienne $(M_\infty^n,g)$ compacte (de même dimension). Alors $M_p$ est difféomorphe à $M_\infty$ pour $p$ assez grand.
\end{theo}

Le but de cet article est de démontrer trois résultats de stabilité nouveaux. Le premier est une version optimisée (et quantitative) du résultat de P.~Petersen \cite{Pe}:

\begin{theo}\label{TheoA}
 Il existe des constantes $\epsilon(n)>0$, $\beta(n)>0$ et $C(n)>0$ telles que si les ${\bf n}$ premières valeurs propres non nulles du laplacien de $M$ sont inférieures à $n+\epsilon(n)$ alors la distance de Gromov-Hausdorff entre $\S^n$ et $M$ est majorée par $C(n)\bigl(\lambda_n-n\bigr)^{\beta(n)}$ et $M$ est donc difféomorphe à $\S^n$.
\end{theo}

Notons que dans notre schéma de preuve, une approximation de Hausdorff à valeur dans $\sn$ est explicitement construite, ce qui rend calculable la valeur $\epsilon(n)$ à partir de laquelle elle existe.

\begin{remark}
  Le théorème A est optimal en ce qui concerne la première conclusion, puisqu'on construit dans cet article une suite de métriques $g_k$ sur $\sn$ telle que $\Ric(g_k)\geq n-1$, $\lambda_i(g_k)\to n$ pour $1\leq i\leq n-1$ et telle que la suite $(\sn,g_k)$ tende (en distance de Gromov-Hausdorff) vers la demi-sphère de dimension $n{-}1$ munie de sa métrique canonique. En revanche, le problème de savoir à partir de quelle valeur $k$, l'hypothèse $\lambda_k(M)\leq n\bigl(1+\epsilon(n)\bigr)$ implique que la variété $M^n$ est difféomorphe à $\sn$ est encore un problème ouvert (il résulte du théorème A et de contre-exemples dûs à M.~Anderson \cite{An2} et Y.~Otsu \cite{Ots} que cette valeur est comprise entre $2$ et $n$).
\end{remark}

 Notons qu'il est démontré dans \cite{Ber} que $M$ possède $k$ petites valeurs propres proches de $n$ si et seulement s'il existe une partie de $M$ Hausdorff proche de $\S^{k-1}$. De ceci et du théorème A découle le:

\begin{corollary}
 Il existe des constantes $\epsilon(n)>0$, $\beta(n)>0$ et $C(n)>0$ telles que, si $M$ contient une partie $A$ qui, munie de la distance induite, est à distance de Gromov-Hausdorff de $(\S^{n-1},can)$ inférieure à $\epsilon\leq\epsilon(n)$, alors $M$ est difféomorphe à $\sn$ et $M$ est à distance de Gromov-Hausdorff de $\sn$ inférieure à $C(n)\epsilon^{\beta(n)}$.
\end{corollary}

\begin{remark}
  Ce corollaire est encore valable s'il existe seulement une partie $A$ de $M$ qui est Gromov-Hausdorff proche du sous ensemble $\{\pm e_1,\ldots,\pm e_n\}$ de $\sn$ (où $(e_1,\ldots,e_n)$ est une base orthonormée de $\sn$).
\end{remark}

Avant d'énoncer les autres résultats principaux de cet article, rappelons le volume d'une variété $M$ de courbure $\Ric\geq n{-}1$ est majoré par $\Vol\sn/\#\pi_1(M)$, l'égalité étant atteinte si et seulement si $\pi_1(M)$ est un sous-groupe (fini) de $O(n+1)$ agissant librement sur $\sn$ et si la variété riemannienne $M$ est l'espace lenticulaire $\sn/\pi_1(M)$. De même, si $M$ est non orientable, alors son volume est majoré par $\Vol(\sn,can)/2$, l'égalité étant atteinte si et seulement si la dimension $n$ est paire et $ M=\mathbb{P}^n\Re$.

Dans le cas où $M$ est de volume presque maximal, nous prouvons les deux résultats suivants (qui améliorent à la fois un résultat de T.~Colding, qui ne traite que le cas où le volume de $M$ est proche de celui de $\sn$ dans \cite{Co1}, et un résultat de T.~Yamaguchi qui suppose de plus que la courbure  sectionnelle de $M$ est minorée par une constante $-K^2$ dans \cite{Yam2}):

\begin{theo}\label{TheoB}
 Il existe des constantes $\epsilon(n)>0$, $\beta(n)>0$ et $C(n)>0$ telles que, si $M$ est non simplement-connexe et vérifie $\Vol M\geq\frac{\Vol\sn}{2}\bigl(1-\epsilon(n)\bigr)$, alors la distance de Gromov-Hausdorff entre $M$ et $\mathbb{P}^n(\Re)$ est majorée par $C(n)\Bigl(\frac{\Vol\sn}{2}-\Vol M\Bigr)^{\beta(n)}$ et $M$ est difféomorphe à $\mathbb{P}^n(\Re)$.
\end{theo}

On remarquera que ce résultat s'applique au cas non-orientable et que, dans ce cas, la conclusion est que $n$ est pair, $M$ est difféomorphe à $\mathbb{P}^n(\Re)$ et Gromov-Hausdorff proche de $\mathbb{P}^n(\Re)$.
Dans le cas où le groupe fondamental est de cardinal quelconque, on a le:

\begin{theo}\label{TheoC}
Soit $k\in\N^*$. Il existe des constantes $\epsilon(n,k)>0$, $\beta(n)>0$ et $C(n)>0$ telles que si $M$ vérifie:
$$\displaylines{\hfill\Vol M\geq\frac{\Vol\sn}{k}\bigl(1-\epsilon(n,k)\bigr)\hfill\mbox{et}\hfill\#\pi_1(M)\geq k,\hfill}$$
alors $\pi_1(M)$ est un sous-groupe (de cardinal $k$) de $O(n+1)$ agissant librement sur $\sn$, la distance de Gromov-Hausdorff entre $ M$ et $\sn/\pi_1(M)$ est majorée par $C(n)\Bigl[\frac{\Vol\sn}{k}-\Vol M\Bigr]^{\beta(n)}$ et $M$ est difféomorphe à l'espace lenticulaire $\sn/\pi_1(M)$.
\end{theo}

En remarquant que tout espace lenticulaire non simplement connexe est de diamètre majoré par $\frac{\pi}{2}$, on obtient, en courbure de Ricci minorée par $n{-}1$, le:

\begin{corollary}
Il existe des constantes $\epsilon(n)>0$ et $\delta(n)>0$  telles que si $M$ vérifie les deux inégalités:
$$\displaylines{\hfill\Vol M\geq\frac{\Vol \sn}{2}\bigl(1-\delta(n)\bigr)\hfill\mbox{et}\hfill\Diam M\geq \frac{\pi}{2}+\epsilon(n),\hfill}$$
alors $M$ est simplement connexe.
\end{corollary}

D'après l'amélioration de l'inégalité de Lichnerowicz sur le $\lambda_1$ donnée par \cite{BBG}, on a un énoncé équivalent où l'hypothèse sur le diamètre est remplacée par l'hypothèse $\lambda_1\leq n+\delta(n)$ (où $\delta(n)$ est une constante strictement positive dont un minorant est donné dans \cite{BBG}).
\bigskip

On obtient aussi l'amélioration suivante des théorèmes de Lichnerowicz et de Myers:
\begin{corollary}
  Il existe des constantes $\alpha(n)>0$ et $\delta(n)>0$  telles que si $M$ est non simplement connexe alors:
$$\displaylines{\lambda_1(M)\geq n\Bigl(\frac{\#\pi_1(M)\Vol M}{\Vol\sn}+\delta(n)\Bigr),\cr
\Diam M\leq\bigl(\frac{\pi}{2}+\delta(n)\bigr)\biggl[1+\frac{1-\frac{\#\pi_1(M)\Vol M}{\Vol\sn}}{\alpha(n)}\biggr].}$$
\end{corollary}

\noindent{\sc Plan de l'article:}
Dans la section 2 nous établissons les estimées analytiques sur les fonctions propres de $M$ nécessaires à nos démonstrations, dans la section 3 nous nous plaçons dans le cas où $M$ admet $n+1$ valeurs propres proches de $n$ et, en suivant \cite{Ga} et \cite{Pe}, nous construisons à l'aide de ces fonctions propres une approximation de Hausdorff $\Phi$ de $M$ dans $\sn$ de degré $\pm1$ (nous prouvons ce dernier fait, qui manquait dans \cite{Pe}; voir l'erratum \cite{Pe2}). Cependant, notre méthode de preuve s'appuiera sur des arguments différents de ceux proposés dans \cite{Ga} et \cite{Pe}. Dans la section 4, nous montrons qu'un sous-groupe fini d'isométries de $M$ admet une action isométrique sur $\sn$ qui rend $\Phi$ équivariante et que, si la première action est libre, alors la seconde l'est aussi. Dans la section 5, nous donnons la preuve des théorèmes B et C. Pour cela, montrons que nos variétés sont Hausdorff proches d'un quotient de la sphère par un groupe fini d'isométries, mais pour pouvoir conclure, il faut toutefois montrer que ce quotient est une variété non singulière. Pour ce point délicat, T.~Yamaguchi utilise dans \cite{Yam2} une hypothèse supplémentaire de minoration de la courbure sectionnelle, qui implique une minoration de la systole de $M$. Nous nous passons de cette hypothèse grace au résultat de la section 4 (quant à la technique de T.~Colding pour traiter le cas $k{=}0$ du théorème \ref{TheoC}, elle utilise de manière essentielle le fait que l'espace modèle est dans ce cas de Radius $\pi$; cette ne peut pas s'adapter trivialement pour démontrer le théorème \ref{TheoC}). En section 6 nous prouvons le théorème A et décrivons en section 7 la famille de contre-exemples qui prouve l'optimalité du théorème A.\\[2mm]
{\sl Remerciements.} Je remercie S.~Gallot pour ses nombreux encouragements et U.~Suter pour les nombreuses conversations concernant le lemme \ref{Suter}.

\section{Valeurs propres proches de $n$ et fonctions propres associées}\label{fibreE}

\noindent{\bf Fibré augmenté de E.~Ruh \cite{Ru}:}
\smallskip

Notons $E\to M$ le fibré vectoriel obtenu comme somme de Withney de $TM$ et d'un fibré trivial en droite (on notera $E=TM\oplus\Re\,e$). Le produit scalaire et la connection linéaire suivants munissent $E$ d'une structure de fibré riemannien:
$$\begin{array}{l}
<X+fe,Y+he>_E=g(X,Y)+fh\\
D^E_Z(X+fe)=D_Z^MX+fZ+\bigl(df(Z)-g(Z,X)\bigr).e
\end{array}$$
où on a noté $g$ la métrique de $M$ et $D^M$ sa connection de Levi-Civita.

\begin{remark}
Le fibré normal de $\sn$ dans $\Re^{n+1}$ est un fibré trivial en droite dont la somme de Withney avec $T\sn$, munie de la métrique décrite plus haut n'est autre que le fibré trivial $E=\sn\times(\rn+1,can)$. Rappelons que les fonctions propres de $\sn$ associées à la valeur propre $n$ sont de la forme $f(x)=<u,x>_{\rn+1}$, où $u$ est un vecteur quelconque de $\rn+1$. Notons $e(x)=x$, alors la correspondance $f\mapsto\nabla f+f.e=u$ donne une identification naturelle entre l'espace des fonctions propres associées à la valeur propre $n$ de $\sn$ et les sections parallèles du fibré $E$. On va montrer dans la suite que cette correspondance se généralise à toutes les variétés de courbure de Ricci presque minorée par $n-1$: le laplacien de $ M$ sur les fonctions aura autant de valeurs propres proches de $n$ que le laplacien brut $\overline{\lap}^E$ sur $E$ aura de valeurs propres proches de $0$ (proposition \ref{klapimpklapbr}) et les sections de $E$ de la forme $\nabla f+fe$, où $f$ est une fonction propre associée à une valeur propre proche de $n$, seront presque "parallèles" (lemme \ref{estimSi}). Inversement les fonctions de la forme $<S,e>_E$, où $S$ sont des sections propres associées à de petites valeurs propres seront presque des fonctions propres associées à des valeurs propres proches de $n$.
\end{remark}

\noindent{\bf Opérateur $\lap_{sph}$:}
\smallskip

Dans la suite, {\sl on note {\boldmath $\pi$} la projection orthogonale de $E$ sur son sous-fibré $TM$ et $A$ la symétrie orthogonale d'hyperplan $TM$. On définit un champ d'endomorphismes symétriques sur $E$ en posant $\Ric'(S)=\Ric_M^{~}\bigl(\pi(S)\bigr)-(n-1)\pi(S)$. On notera par la suite $\lap_{sph}$ l'opérateur $\overline{\lap}^E+\Ric'$}.

\noindent{\bf Sections $S_f=\nabla f+f.e$:}
\smallskip

Le lemme suivant sera fondamental pour relier les fonctions propres de $ M$ aux sections propres de $\overline{\lap}^E$~:

\begin{lemma}\label{laplapsph}
  Soit $f:M\to \Re$ telle que $\lap f=\lambda f$. On pose $S_f=\nabla f+f.e$, on a alors la relation $\lap_{sph}(S_f)=(\lambda-n)A(S_f)$.
\end{lemma}

\begin{proof}
  Un calcul direct donne:\\
$$\overline{\lap}^E(S_f)(x)=\overline{\lap}_M\nabla f-\nabla f+nf.e-\lap f.e,$$
où $\overline{\lap}_M$ est le laplacien brut de $TM$. Enfin, l'opérateur $\overline{\lap}_M+\Ric_M$ n'est autre que l'opérateur de Hodge $\lap_H$ sur les champs de vecteurs, qui commute avec $\nabla$, d'où la relation annoncée.
\end{proof}

\noindent{\bf Correspondance fonctions/sections propres}
\smallskip

Nous aurons besoin d'une minoration de la première valeur propre $\lambda_1^1$ du laplacien de Hodge sur les 1-formes de $M$. L'inégalité suivante découle d'une variante de la méthode de Lichnerowicz. Nous en donnons une preuve succinte:
\smallskip

\begin{lemma}\label{L1}
  Soit $ M$ vérifiant $\Ric\geq n-1$, alors $\lambda_1^1(M)\geq n$ (resp. la première valeur propre du laplacien de Hodge, restreint aux 1-formes co-fermées est minorée par $2(n-1)$).
\end{lemma}

\begin{proof}
  Soit $\alpha\in\Lambda^1(M)$. En décomposant orthogonalement $D\alpha$ en partie antisymétrique $\frac{d\alpha}{2}$, partie symétrique sans trace et partie scalaire $-\frac{\delta\alpha}{n}g$, on obtient $|D\alpha|^2\geq\frac{|d\alpha|^2}{2}+\frac{(\delta\alpha)^2}{n}$.
En intégrant la formule de Bochner on obtient:
$$(\lap\alpha,\alpha)=\|D\alpha\|_2^2+\int_M\Ric(\alpha,\alpha)\geq\frac{\|d\alpha\|_2^2}{2}+\frac{\|\delta\alpha\|_2^2}{n}+(n-1)\|\alpha\|_2^2$$
Si $d\alpha=0$, alors $\|\delta\alpha\|_2^2=(\lap\alpha,\alpha)$ et donc $(\lap\alpha,\alpha)\geq n\|\alpha\|^2$. Si $\delta\alpha=0$, alors  $\|d\alpha\|_2^2=(\lap\alpha,\alpha)$ et donc $(\lap\alpha,\alpha)\geq 2(n-1)\|\alpha\|^2$.
\end{proof}

Nous en déduisons que $\lap$, $\overline{\lap}^E$ et $\lap_{sph}$ ont le même nombre de petites valeurs propres, c'est la:
\smallskip

\begin{proposition}\label{klapimpklapbr}
   Il existe des constantes $\alpha'(n)>\alpha(n)>0$ (explicitement calculables) telles que, si $ M$ est une variété complète vérifiant $\Ric\geq n{-}1$, alors:
\medskip

\noindent $(i)$ $\lambda_{n+2}(M)\geq n+\alpha'(n)$ et $\lambda_{n+2}(\lap_{sph})\geq\lambda_{n+2}(\overline{\lap}^E)\geq \alpha'(n)$.
\medskip

\noindent $(ii)$ Pour tout entier $k$, $0\leq\lambda_{k}(\overline{\lap}^E)\leq\lambda_k(\lap_{sph})\leq\bigl(\lambda_k(M)-n\bigr).$

Réciproquement, si $\lambda_{k}(\overline{\lap}^E)\leq\alpha(n)$ alors $$\lambda_{k}(M)\leq n+200n^2\sqrt{\lambda_{k}(\overline{\lap}^E)}.$$
\end{proposition}

\begin{proof}
La première série d'inégalités annoncée en $(ii)$ découle directement de la positivité du potentiel $\Ric'$, du lemme \ref{laplapsph} et du principe du min-max.

Pour finir de démontrer $(ii)$, notons $(S_i)_{1\leq i\leq k}$ une famille $L^2$-ortho--normée de sections propres (associées aux $k$ premières valeurs propres) de $\overline{\lap}^E$. Notons ${\cal E}$ l'espace vectoriel engendré par les sections $(S_i)_{1\leq i\leq k}$, muni du produit scalaire $L^2$ de $E$, et ${\cal F}=\Psi({\cal E})$, où $\Psi(S)=<S,e>_E$ (notons que ${\cal F}$ est engendré par les fonctions $f_i=<S_i,e>_E$). Nous allons conclure en appliquant le principe du min-max à ${\cal F}$:

$\bullet$ Les fonctions de ${\cal F}$ sont d'intégrale nulle car $\overline{\lap}^E$ est auto-adjoint, $\overline{\lap}^Ee=n.e$ (lemme \ref{laplapsph}, avec $f=1$), $\overline{\lap}^E(S_i)=\lambda_iS_i$ et $\lambda_i\neq n$, d'où $\frac{1}{\Vol M}\int_Mf_i=\frac{1}{\Vol M}\int_M<S_i,e>_E^{~}=0$.

$\bullet$ $\Psi$ est une application injective, et donc ${\cal F}$ est un espace de fonctions de dimension $k$: si $S\in{\rm Ker}\Psi\cap{\cal E}$, alors $S=X\in\Gamma(M)$ et  $\|D^ES\|_2^2\leq\alpha(n)\|S\|_2^2<\|S\|_2^2$. Or $D^E_YX=D^M_YX-g(X,Y)e$, et donc $\|D^ES\|_2^2=\|D^MX\|_2^2+\|X\|_2^2\geq\|X\|_2^2=\|S\|_2^2.$ On en déduit que $S=0$.

$\bullet$ Il ne reste plus qu'à montrer que le quotient de Rayleigh des fonctions de ${\cal F}$ est presque plus petit que $n$. Soit donc $f\in{\cal F}\setminus\{0\}$. Il existe donc $\alpha\in\Lambda^1(M)$ tel que $S=\alpha^\#+f.e\in{\cal E}$, ce qui donne:
$$\displaylines{\|D^M\alpha^\#+fId_{TM}-(df-\alpha).e\|_2^2=\|D^ES\|^2_2\hfill\cr
\hfill\leq\lambda_k(\overline{\lap}^E)\|S\|^2_2=\lambda_k(\overline{\lap}^E)(\|\alpha\|_2^2+\|f\|_2^2),}$$
Dont on déduit les deux inégalités suivantes:
$$\|D^M\alpha+fg\|_2^2\leq\lambda_k(\overline{\lap}^E)(\|\alpha\|_2^2+\|f\|_2^2)$$
$$\|df-\alpha\|_2^2\leq\lambda_k(\overline{\lap}^E)(\|\alpha\|_2^2+\|f\|_2^2)$$
La deuxième de ces inégalités donne l'inégalité:
$$\Bigl|\frac{\|df\|_2}{\|f\|_2}-\frac{\|\alpha\|_2}{\|f\|_2}\Bigr|\leq\sqrt{\lambda_k(\overline{\lap}^E)}\Bigl(1+\frac{\|\alpha\|_2}{\|f\|_2}\Bigr).$$
On est donc ramené à montrer que le rapport $\frac{\|\alpha\|_2^2}{\|f\|_2^2}$ est presque majoré par $n$. Or, en procédant comme dans la démonstration du lemme \ref{L1}:
$$\frac{\|d\alpha\|_2^2}{2}+\|f-\frac{\delta\alpha}{n}\|_2^2\leq\|D^M\alpha+fg\|_2^2\leq\lambda_k(\overline{\lap}^E)(\|\alpha\|_2^2+\|f\|_2^2).$$
On en déduit:
$$\displaylines{(\lap\alpha,\alpha)=\|d\alpha\|_2^2+\|\delta\alpha\|_2^2\hfill\cr
\hfill\leq5n^2\sqrt{\lambda_k(\overline{\lap}^E)}(\|\alpha\|_2^2+\|f\|_2^2)+\Bigr(1+\sqrt{\lambda_k(\overline{\lap}^E)}\Bigl)n^2\|f\|_2^2.}$$
Enfin, d'après le lemme \ref{L1}, on a $(\lap\alpha,\alpha)\geq n\|\alpha\|_2^2$. En combinant les deux dernières inégalités, on obtient $\|\alpha\|_2^2\leq\Bigl(1+12n\sqrt{\lambda_k(\overline{\lap}^E)}\Bigr)n\|f\|_2^2,$
dès que $\alpha(n)\leq\frac{1}{36n^4}$, ce qui permet de conclure.

$(i)$ sera démontré pour l'opérateur $\overline{\lap}^E$ dans la section suivante. Pour les autres opérateurs, cela découle alors de $(ii)$.
\end{proof}

\begin{remark}
  On pourrait renforcer l'analogie avec le cas de $\sn$ décrit précédemment en montrant que, si $f$ est une combinaison linéaire des fonctions propres de $M$ associées à des valeurs propres proches de $n$, alors $S_f$ est $L^2$-proche d'une combinaison linéaire des sections propres de $\overline{\lap}^E$ associées à de petites valeurs propres (la réciproque est aussi vraie en remplaçant $f\mapsto S_f$ par $S\mapsto<S,e>_E$).
\end{remark}

\noindent{\bf Estimées sur les fonctions propres}
\smallskip

Soit $(\sqrt{n+1}.f_i)_{1\leq i\leq k}$ une famille $L^2$-orthonormée de fonctions pro--pres de $M$ associées à des valeurs propres $0<\lambda_1\leq\ldots\leq\lambda_k\leq n+\epsilon$. On lui associe la famille $L^2$-orthogonale $(S_i)_{1\leq i\leq k}$ de sections du fibré $E$ définies par $S_i=\nabla f_i+f_ie$. 
 On obtient alors les estimées analytiques et géométriques suivantes~:
\smallskip

\begin{lemma}\label{estimSi}
  Il existe des constantes $\alpha(n)>0$ et $C(n)>0$ (explicitement calculables) telles que si $M$ vérifie $\lambda_{k}\leq n+\epsilon$ (avec $\epsilon\leq\alpha(n)$), alors~:\\
$(i)\hskip1cm\|\sum_{i=1}^{k}\alpha_iS_i\|_\infty\leq\Bigl(1+C(n)\sqrt{\epsilon}\Bigr)\|\sum_{i=1}^{k}\alpha_iS_i\|_2,$\\
pour tout $(\alpha_i)\in\Re^{k}$.\\
$(ii)$ Il existe un sous-ensemble $M_{\epsilon}$ de $M$ tel que~:\\
~\hskip5mm $\Vol M_{\epsilon}\geq\bigl(1-C(n)\epsilon^\frac{1}{4}\bigr)\Vol M$\\
\hskip5mm$|<S_i(x),S_j(x)>_E-\delta_{ij}|\leq C(n)\epsilon^\frac{1}{4}$ pour tout $x\in M_{\epsilon}$ et tout couple $(i,j)$.
\end{lemma}
\begin{remark}
   Si la famille $(S_i)_{1\leq i\leq k}$ est une famille $L^2$-orthonormée de sections propres d'un opérateur $\overline{\lap}^E+V$ à potentiel $V$ positif associées à des valeurs propres $\lambda_1(\overline{\lap}^E+V)\leq\ldots\leq\lambda_k(\overline{\lap}^E+V)\leq\epsilon$ (avec $\epsilon\leq\alpha(n)$), alors les propriétés $(i)$ et $(ii)$ sont encore valables (la preuve qui suit s'adapte facilement).
\end{remark}

\begin{proof}
D'après le lemme \ref{laplapsph}, la famille de sections  $(\sqrt{\frac{n+1}{\lambda_i+1}}S_i)$ est $L^2$-orthonormée et vérifie $\lap_{sph}S_i=(\lambda_i-n)A(S_i)$ pour tout indice $i\leq k$. Soient $E_k$ l'espace vectoriel engendré par $(S_i)_{1\leq i\leq k}$ et $\dis A_p=\supp_{S\in E_k\setminus\{0\}}\frac{\|S\|_p}{\|S\|_2}$, où $p\in]1,+\infty]$. 

Soit $S\in E_k\setminus\{0\}$ et posons $u=\sqrt{|S|^2+\epsilon^2}$ pour $\epsilon>0$. D'après l'inégalité de Kato, on a $u\lap u\leq\,<\overline{\lap}^ES,S>_E\leq|\lap_{sph}S|u$. On en déduit que, pour tout réel $p>1/2$~:
$$\|d(u^p)\|_2^2=\frac{p^2}{2p-1}\int_M(u\lap u)u^{2p-2}\leq\frac{p^2}{2p-1}\left[\|\lap_{sph}S\|_{2p}\|u\|^{2p-1}_{2p}\right]$$
En appliquant l'inégalité de Sobolev $\|f\|_\frac{2n}{n-2}^2\leq C(n)\|df\|_2^2+\|f\|_2^2$ donnée par \cite{Ili} à la fonction $u^p$ et en faisant tendre $\epsilon$ vers 0, on obtient~:
$$
\| S\|_{\frac{2pn}{n-2}}^p\leq\frac{C(n)p}{\sqrt{2p-1}}\sqrt{\|\lap_{sph}S\|_{2p}\| S\|_{2p}^{2p-1}}+\| S\|_{2p}^p$$
(cette inégalité reste valable en dimension $2$ en remplaçant $n$ par $4$).
Or $E_k$ est un espace stable par $A\lap_{sph}$ et $A$ est une isométrie, on a donc~:
$$\|\lap_{sph}S\|_{2p}=\|A\lap_{sph}S\|_{2p}\leq A_{2p}\|A\lap_{sph}S\|_2\leq A_{2p}(\lambda_k-n).\|S\|_2.$$
On en déduit l'inégalité $A_{\frac{2pn}{n-2}}\leq\left(1+\frac{C(n)p}{\sqrt{2p-1}}(\lambda_k-n)^\frac{1}{2}\right)^{1/p}A_{2p}$, et donc $A_\infty\leq\prod_{j=0}^{\infty}\left(1+\frac{C(n)\nu^j}{\sqrt{2\nu^j-1}}(\lambda_k-n)^\frac{1}{2}\right)^{\frac{1}{\nu^j}}$ en posant $\nu=\frac{n}{(n-2)}$.
Enfin, en utilisant la concavité de la fonction $\log$, on obtient $A_\infty\leq\bigl(1+C(n)\sqrt{\epsilon}\bigr)$.

Nous allons maintenant prouver que l'inégalité $(i)$ implique la propriété $(ii)$. Posons $M_\epsilon$ l'ensemble des points $x\in M$ où:
$$\displaylines{\hfill|S_i(x)+S_j(x)|_E^2\geq2(1-\epsilon^\frac{1}{4}),\hfill|S_i(x)-S_j(x)|_E^2\geq2(1-\epsilon^\frac{1}{4}),\hfill\cr
\hfill|S_i(x)|_E^2\geq1-\epsilon^\frac{1}{4}\hfill\forall i<j\hfill}$$
D'après $(i)$ (et en posant $h=\sum_{i<j}|S_i+S_j|_E^2+|S_i-S_j|_E^2+\sum_i2|S_i|_E^2$)~:
$$\displaylines{
2k^2=\frac{1}{\Vol M}\int_Mh=\frac{1}{\Vol M}\int_{M_\epsilon}h+\frac{1}{\Vol M}\int_{M\setminus M_\epsilon}h\hfill\cr
\leq\frac{\Vol M_\epsilon}{\Vol M}\Bigl(\sum_{i<j}\|S_i+S_j\|_\infty^2+\|S_i-S_j\|_\infty^2+\sum_i2\|S_i\|_\infty^2\Bigr)\hfill\cr
\hfill+(k^2-1)\Bigl(1-\frac{\Vol M_\epsilon}{\Vol M}\Bigr)\max_{1\leq i<j\leq k}\Bigl[\|S_i+S_j\|_\infty^2,\|S_i-S_j\|_\infty^2,2\|S_i\|_\infty^2\Bigr]\hfill\cr
+2\Bigl(1-\frac{\Vol M_\epsilon}{\Vol M}\Bigr)(1-\epsilon^\frac{1}{4})\cr
\leq2k^2\bigl(1+C(n)\sqrt{\epsilon}\bigr)\frac{\Vol M_\epsilon}{\Vol M}\hfill\cr
\hfill+\Bigl[2(k^2-1)\bigl(1+C(n)\sqrt{\epsilon}\bigr)+2(1-\epsilon^\frac{1}{4})\Bigr]\left(1-\frac{\Vol M_\epsilon}{\Vol M}\right)
}$$
On en déduit $\frac{\Vol M_\epsilon}{\Vol M}\geq1-\frac{k^2C(n)\epsilon^\frac{1}{4}}{1+C(n)\epsilon^\frac{1}{4}}\geq1-k^2C(n)\epsilon^\frac{1}{4}$, pour $\epsilon\leq\alpha(k,n)$ assez petit.
D'après $(i)$, on a pour tout $x\in M_\epsilon$ et tout $i\neq j$~:
\begin{eqnarray*}
  \bigl|<S_i(x),S_j(x)>_E\bigr|&=&\Bigl|\frac{|S_i(x)+S_j(x)|^2-|S_i(x)-S_j(x)|^2}{4}\Bigr|\\
&\leq&\frac{2(1+C(n)\sqrt{\epsilon})-2(1-\epsilon^\frac{1}{4})}{4}\leq C'(n)\epsilon^\frac{1}{4}
\end{eqnarray*}
De même, $\bigl|\|S_i(x)\|_E^2-1\bigr|\leq C'(n)\epsilon^\frac{1}{4}$. Supposons $k>n+1$: en appliquant ce qui précède à la famille $(S_i)_{1\leq i\leq n+2}$, on obtient que $|<S_i,S_j>-\delta_{ij}|\leq C(n)\epsilon^\frac{1}{4}$ sur $M_\epsilon$. Mais alors, il existe au moins un point $x$ de $M$ où la famille $(S_i(x))_{1\leq i\leq n+2}$ est de rang $n+2$, ce qui est impossible car $E$ est de rang $n+1$. On en déduit d'abord que $k\leq n+1$ (et donc le $(i)$ de la proposition \ref{klapimpklapbr}), puis la propriété $(ii)$ annoncée.\end{proof}

\section{Variétés admettant $n+1$ petites valeurs propres}\label{Ln+1}

Dans cette section nous étudions les variétés complètes vérifiant $\Ric\geq(n{-}1)$ et $\lambda_{n+1}\leq n+\alpha(n)$. Fixons une famille $(\sqrt{n+1}.f_i)_{1\leq i\leq n+1}^{~}$ $L^2$-orthonormée de fonctions propres associées aux valeurs propres $n\leq\lambda_1\leq\ldots\leq\lambda_{n+1}\leq n+\epsilon$ et considérons l'application:
\begin{eqnarray}\label{Fi}
\Phi~: M&\to&\S^n\hookrightarrow\Re^{n+1}\nonumber\\
x&\mapsto&\frac{1}{\bigl(\sum_jf_j(x)^2\bigr)^{1/2}}.\bigl(f_1(x),\ldots,f_{n+1}(x)\bigr)
  \end{eqnarray}
Nous démontrons dans cette section la:
\smallskip

\begin{proposition}\label{Fiprop}
   Il existe des constantes $\alpha(n)>0$ et $C(n)>0$ (explicitement calculables) telles que si $ M$ est une variété complète qui vérifie $\Ric\geq(n-1)$ et $\lambda_{n+1}\leq n+\epsilon$ (pour $\epsilon\leq\alpha(n)$), alors~:
$$\Vol M\geq\bigl(1-C(n)\epsilon^\frac{1}{2(n+1)}\bigr)\Vol\sn.$$
et l'application $\Phi$ est une approximation de Hausdorff à valeur dans $\sn$ de degré $\pm 1$.
Plus précisement, pour tout couple de points $(x,y)$ de $M$, on a:
 $$\left|d_{\sn}(\Phi(x),\Phi(y))-d_M(x,y)\right|\leq C(n)\epsilon^\frac{1}{384(n+1)^3}.$$
\end{proposition}

\begin{remark} 
Notre schéma de preuve est une adaptation de celui de P.~Petersen dans \cite{Pe} (voir aussi \cite{Au1}). Toutefois, nous corrigeons l'erreur faite par P.~Petersen sur le calcul du degré de $\Phi$ (voir l'erratum \cite{Pe2}) ce qui est indispensable pour notre démonstration des résultats annoncés en introduction. Pour la démonstration du fait que $\Phi$ est une approximation de Hausdorff, nous renvoyons à l'article \cite{Pe} de P.~Petersen (noter que d'après ce qui suit, $\Phi$ est surjective). On pourra regarder aussi \cite{Au2}, où la même propriété est démontrée sous des hypothèses de pincement intégral de la courbure de Ricci passant sous $n{-}1$ et la forme du majorant de la distance de Gromov-Hausdorff est explicitement calculée.
\end{remark}

\noindent{\bf $\Phi$ est bien définie sous nos hypothèses~:}
\smallskip

Cela découle du:
\smallskip

\begin{lemma}\label{1.13}
 Il existe des constantes $\alpha(n)>0$ et $C(n)>0$ (explicitement calculables) telles que si $\lambda_{n+1}\leq n+\epsilon$ (avec $\epsilon\leq\alpha(n)$), on a~:
$\|\sum_{i=1}^{n+1}f^2_i~-1\|_\infty\leq C(n)\epsilon^{\frac{1}{2(n+1)}}.$\end{lemma}

\begin{proof}
Soit $x_0$ un point quelconque de $M$. En appliquant le lemme \ref{estimSi} avec $\alpha_i=f_i(x_0)$, on obtient:
$$\displaylines{\bigl(\sum_{i=1}^{n+1}f_i(x_0)f_i\bigr)^2+\bigl(\sum_{i=1}^{n+1}f_i(x_0)\nabla f_i\bigr)^2\leq(1+C(n)\sqrt{\epsilon})^2\|\summ_{i=1}^{n+1}f_i(x_0)S_i\|_2^2\cr
\hfill\leq(1+C(n)\sqrt{\epsilon})\sum_{i=1}^{n+1}f_i(x_0)^2.}$$
On en déduit que la fonction $h=\summ_{i=1}^{n+1}f_i^2$ vérifie les inégalités $\|h\|_1=1$, $\|h\|_\infty\leq\bigl(1+C'(n)\sqrt{\epsilon}\bigr)$ et $\|dh\|_\infty\leq C'(n)$. Le résultat annoncé découle alors du lemme suivant~:

\begin{lemma}
  Soient $(M^n,g)$ compacte et vérifiant $\Ric\geq0$ et $h:M\to\Re$ une fonction lipschitzienne, alors:
$$|h|\geq\|h\|_\infty-2\bigl(\Diam M\|dh\|_{\infty}\bigr)^\frac{n}{n+1}\bigl(\|h\|_\infty-\|h\|_1\bigr)^\frac{1}{n+1}.$$
\end{lemma}

\begin{proof}
On peut supposer $\|h\|_1>0$. Soit $x\in M$:
$$\displaylines{\|h\|_1=\frac{1}{\Vol M}\int_{B(x,\eta)}|h|+\frac{1}{\Vol M}\int_{M\setminus B(x,\eta)}|h|\hfill\cr
\hfill\leq\frac{\Vol B(x,\eta)}{\Vol M}\bigl(|h(x)|+\|df\|_\infty\eta\bigr)+\left(1-\frac{\Vol B(x,\eta)}{\Vol M}\right)\|h\|_\infty.\hfill}$$
On en déduit que $\|h\|_\infty-\|h\|_1\geq\frac{\Vol B(x,\eta)}{\Vol M}\bigl(\|h\|_\infty-|h(x)|-\|df\|_\infty\eta\bigr)$. Posons $\eta=\frac{\|h\|_\infty-|h(x)|}{2\|df\|_\infty}>0$. Le théorème de Bishop-Gromov, donne alors:
$$\|h\|_\infty-\|h\|_1\geq\frac{\bigl(\|h\|_\infty-|h(x)|\bigr)^{n+1}}{2^{n+1}\|df\|_\infty^n(\Diam M)^n}.$$
\end{proof}
\end{proof}

\noindent{\bf Calcul de $d\Phi$~:}

Soit $x$ un point de $M$, alors $d_x\Phi$ est une application de $T_xM$ dans $T_{\Phi(x)}\S^n\subset\Re^{n+1}$. On note $^td_x\Phi$ l'application transposée de $d_x\Phi$ vu comme une application de $(T_xM,g_x)$ dans $\Re^{n{+}1}$ muni de son produit scalaire canonique, qu'on notera $<.,.>^{~}_{\Re^{n+1}}$. Si $(\varepsilon_i)_{1\leq i\leq n+1}$ est la base canonique de $\Re^{n+1}$ alors $^td_x\Phi(\varepsilon_i)=\frac{\nabla f_i-\Phi_i\sum_{k=1}^{n+1}\Phi_k.\nabla f_k}{(\sum_j f_j^2)^{\frac{1}{2}}}=\frac{S_i-\Phi_i\sum_{j=1}^{n+1}\Phi_jS_j}{(\sum_j f_j^2)^{\frac{1}{2}}}$, où on a posé $\Phi_i=\frac{f_i}{\bigl(\sum_jf_j^2\bigl)^\frac{1}{2}}$ (où $TM$ est vu comme un sous-fibré de $E$). En particulier, pour tout vecteur $v$ de $T_{\Phi(x)}\S^n(=\Phi(x)^{\bot})$, on a $^td_x\Phi(v)=\frac{\sum_{i=1}^{n+1}v_iS_i}{(\sum f_j^2)^{\frac{1}{2}}}$.

\noindent{\bf $\Phi$ est presque contractante:}
\smallskip

Les lemmes \ref{estimSi} et \ref{1.13} nous donnent alors le:
\smallskip

\begin{lemma}\label{Ficontractante}
 Si $ M$ vérifie $\Ric\geq n-1$ et si $\lambda_{n+1}\leq n+\epsilon$ (avec $\epsilon\leq\alpha(n)$), alors $\|d\Phi\|_\infty\leq1+C(n)\epsilon^{\frac{1}{2(n+1)}}$.
\end{lemma}

\noindent{\bf Si $M$ est orientable, alors ${\rm deg}\,\Phi=\pm1$}
\smallskip

Soit $x$ un point de $M$ et $(X_i)$ une base orthonormée directe de $T_xM$. Pour tout $x\in M$, on munit $E_x$ de l'orientation induite par $TM$ (i.e. telle que $(\tilde{X_i})_{1\leq i\leq n+1}=(X_1,\cdots,X_n,e)$ soit une base directe de $E_x$). On note $L_x:\Re^{n+1}\to E_x$ l'application définie par $L_x(v)=\sum_{i=1}^{n+1}v_iS_i(x)$ et $h(x)=\det L_x$ (le déterminant étant calculé relativement à des bases orthonormées directes des espaces de départ et d'arrivée). En calculant $h$ dans une base orthonormée directe de $\Re^{n+1}$ de la forme $(e_1,\cdots,e_n,\Phi(x))$, on obtient:
$$h(x)=\Bigl(\sum_jf_j^2(x)\Bigr)^\frac{n+1}{2}\det d_x\Phi.$$
Commençons par estimer $\|h\|_2^{~}$. On a $\|h\|^2_\infty\leq\bigl(1+C(n)\sqrt{\epsilon}\bigr)$ car $h^2(x)\leq|L_x|^{2n}$ et d'après le lemme \ref{estimSi} $(i)$, on a: 
$$|L_x(v)|^2_E\leq\bigl(1+C(n)\sqrt{\epsilon}\bigr)\|\sum_{i=1}^{n+1}v_iS_i\|_2^2\leq\bigl(1+C(n)\sqrt{\epsilon}\bigr)\|v\|^2_{\Re^{n+1}},$$
pour tout $v\in\Re^{n+1}$. Par ailleurs, pour tout couple $(u,v)$ de vecteurs de $\Re^{n+1}$, on a~:
$$\displaylines{\Bigl|<^tL_x\circ L_x(u),v>_{\rn+1}^{~}{-}<u,v>_{\rn+1}^{~}\Bigr|{=}\Bigl|\sum_{ij}u_iv_j\bigl(<S_i,S_j>_E-\delta_{ij}\bigr)\Bigr|\hfill\cr
\hfill\leq\max_{ij}\bigl|<S_i,S_j>_E-\delta_{ij}\bigr|\|u\|_{\rn+1}^{~}\|v\|_{\rn+1}^{~}}$$
et donc, d'après le lemme \ref{estimSi} $(ii)$, $\bigl\| ^tL_x\circ L_x-Id_{T_{\Phi(x)}\S^n}\bigr\|\leq C(n)\epsilon^\frac{1}{4}$, pour tout point $x$ du sous-ensemble $M_\epsilon$, i.e. $|h^2(x)-1|\leq C(p,n)\epsilon^\frac{1}{4}$ sur $M_\epsilon$. On obtient alors l'estimée~:
\begin{eqnarray*}
  \Bigl|\frac{1}{\Vol M}\int_Mh^2-1\Bigr|&\leq&\Bigl|\frac{1}{\Vol M}\int_{M_\epsilon}(h^2-1)\Bigr|\\
&&\hskip1.5cm+\Bigl(1-\frac{\Vol M_\epsilon}{\Vol M}\Bigr)\max\bigl(1,\|h^2\|_\infty-1|\bigr)\\
&\leq&C(n)\epsilon^\frac{1}{4}
\end{eqnarray*}

Pour obtenir le degré de $\Phi$, c'est l'intégrale $\frac{1}{\Vol M}\int_Mh$ que l'on doit estimer.Nous allons utiliser l'inégalité de Poincaré pour montrer que la moyenne de $h$ est proche de sa norme $L^2$. Pour cela, il faut montrer que la norme $L^2$ de $dh$ est petite:

 En prenant la base canonique de $\Re^{n+1}$ comme base orthonormée au départ, et une base orthonormée directe quelconque $(\tilde{X}_i(x))$ de $E_x$ à l'arrivée, on obtient $h(x)=det\bigl(<\tilde{X}_i,S_j>_E\bigr)$. Soit $x$ un point donné de $M$, $X\in T_xM$ et $\gamma_X$ la géodésique passant par $x$ avec le vecteur vitesse $X$. On note $(\tilde{X}_i)$ le repère transporté parallèlement le long de $\gamma_X$ de $(\tilde{X}_i(x))$. Alors $d_xh(X)=X.\Bigl(det\bigl(<\tilde{X}_i,S_j>_E^{~}\bigr)_{ij}\Bigr)$ est égal à:
$$
\sum_{j=1}^{n+1}det
\begin{pmatrix}
  {<}\tilde{X}_1,S_1{>}&\hskip-3mm\ldots&\hskip-3mm {<}\tilde{X}_1,D^E_XS_j{>}\hskip-3mm&\hskip-3mm\ldots&\hskip-3mm{<}\tilde{X}_1,S_{n+1}{>}\\
\Biggl|&&\Biggl|&&\Biggl|\\
{<}\tilde{X}_{n+1},S_1{>}&\hskip-3mm\ldots&\hskip-3mm {<}\tilde{X}_{n+1},D^E_XS_j{>}&\hskip-3mm\ldots&\hskip-3mm{<}\tilde{X}_{n+1},S_{n+1}{>}
\end{pmatrix}
$$
D'où $\bigl|d_xh(X)\bigr|\leq C(n)\max_i\|S_i\|_\infty^n\max_i|D^E_XS_i(x)|_E$. Le lemme \ref{estimSi} $(i)$, donne donc:
$$
  \|dh\|_2^2=\frac{1}{\Vol M}\int_M|dh|^2\leq C(n)\max_i\frac{1}{\Vol M}\int|D^ES_i|^2(x)\leq C(n)\epsilon,$$
la dernière inégalité découle de la preuve du lemme \ref{klapimpklapbr} $(ii)$. 

Or $\lambda_1\geq n$ et donc $0\leq\|h\|_2^2-\bigl(\frac{1}{\Vol M}\int_Mh\bigr)^2\leq C(n)\epsilon$. On déduit donc de l'estimée précédente sur $\|h\|_2^2$ que~:
$$\displaylines{1-C(n)\epsilon^\frac{1}{4}\leq\Bigl(\frac{1}{\Vol M}\int_Mh\Bigr)^2\leq1+ C(n)\epsilon^\frac{1}{4}}$$
Pour conclure, on a ${\rm deg}\Phi\Vol\S^n=\int_M\det\,d\Phi=\int_M\bigl(\sum_kf^2_k\bigr)^{-\frac{n+1}{2}}h$. Le lemme \ref{1.13} et la majoration de $\|h\|_\infty$ donnent alors $\Bigl||{\rm deg}\,\Phi|\frac{\Vol\S^n}{\Vol M}-1\Bigr|\leq C(n)\epsilon^\frac{1}{2(n+1)}.$
Comme ${\rm deg}\,\Phi$ est un entier et que $\Vol M\leq\Vol\sn,$ on obtient ${\rm deg}\Phi=\pm1$ et $\Vol M\geq\Vol\S^n\bigl(1-C(n)\epsilon^\frac{1}{2(n+1)}\bigr)$.\medskip

\noindent{\bf Si $\lambda_{n+1}(M)\leq n+\epsilon(n)$ alors $M$ est orientable:}
\smallskip

Si $M^n$ est non-orientable, on construit l'application $\Phi$ à partir des fonctions propres $(f_i)_{1\leq i\leq n+1}$ et on note $\pi:\widetilde{M}\to M$ le revêtement riemannien orientable de $M$. Soit $\widetilde{\Phi}=\Phi\circ\pi:\widetilde{M}\to\S^n$. On a $\deg_2\widetilde{\Phi}=\deg_2\Phi\deg_2\pi=0$ (où $\deg_2$ désigne le degré modulo 2), et donc $\widetilde{\Phi}$ est de degré orientable pair. Cependant, $(\widetilde{M}^n,\tilde{g})$ vérifie aussi $\Ric\geq(n{-}1)$ et les fonctions $\tilde{f}_i=f_i\circ\pi$ sont des fonctions propres  de $(\widetilde{M}^n,\tilde{g})$ associées aux mêmes valeurs propres $\lambda_i$. On en déduit que l'application $\widetilde{\Phi}$ n'est autre que celle étudiée en {\bf e)} associée à la variété $(\widetilde{M},\tilde{g})$ et aux fonctions propres $\tilde{f}_i$. $\widetilde{\Phi}$ doit donc être de degré $\pm 1$, ce qui est contradictoire.

\section{Équivariance de $\Phi$}
\label{sec:Equiv}
Nous allons exhiber, sous les hypothèses de la section \ref{Ln+1}, un morphisme du groupe ${\rm Isom} M$ des isométries de $ M$ dans le groupe $O_{n+1}(\Re)$ des isométries de $\sn$ pour lequel l'application $\Phi$ est équivariante:

Soit $f$ une fonction propre ($\lap f=\lambda f$) et $\sigma$ une isométrie de la variété $ M$. Alors $\lap(f\circ\sigma)=\lambda(f\circ\sigma)$. En particulier, si on suppose que $(\sqrt{n+1}.f_i)_{1\leq i\leq k}$ est une base de l'espace de toutes les fonctions propres associées aux valeurs propres de $ M$ inférieures à un réel $\alpha$ donné, alors il existe une matrice $A^\sigma$ de ${\cal M}_k(\Re)$ telle que pour tout indice $i$, on ait $f_i\circ\sigma=\sum_jA^\sigma_{ij}f_j$.
De plus, comme $\sigma$ préserve le volume riemannien, on a $\frac{\delta_{ij}}{n+1}=\frac{1}{\Vol M}\int_Mf_i\circ\sigma.f_j\circ\sigma dv_g=\sum_{kl}A^\sigma_{ik}A^\sigma_{jl}\frac{\delta_{kl}}{n+1}$.
On en déduit que $\sigma\mapsto A^\sigma$ est un morphisme de groupe de ${\rm Isom M}$ à valeur dans $O_k(\Re)$.

Dans le cas particulier où $ M$ vérifie $\Ric\geq(n-1)$ et $\lambda_{n+1}\leq n+\epsilon$, le choix d'une famille $(\sqrt{n+1}.f_i)_{1\leq i\leq n+1}$ $L^2$-orthonormée de fonctions propres associées aux $n+1$ premières valeurs propres non nulles de $M$ induit un morphisme de ${\rm Isom} M$ à valeurs dans $O(n+1)$ (d'après le lemme \ref{klapimpklapbr} $(i)$) tel que $\Phi\bigl(\sigma(x)\bigr)=A^\sigma\bigl(\Phi(x)\bigr)$ pour toute isométrie $\sigma$. L'outil central de la démonstration des théorèmes B et C est le lemme suivant:

\begin{lemma}\label{quotient}
  Soient $ M$ une variété complète vérifiant $\Ric\geq n{-}1$ et $\lambda_{n+1}\leq n+\epsilon$ (pour $\epsilon\leq\alpha(n)$), et $\Phi:M\to\sn$ l'application définie dans la section \ref{Ln+1}. On se donne de plus $G$ un groupe fini d'isométries agissant librement sur $ M$. Le morphisme $\sigma\mapsto A^\sigma$ est alors injectif sur $G$ et $A(G)$ est un sous groupe de $O_{n+1}(\Re)$ agissant librement sur $\sn$. Enfin, $\Phi$ passe au quotient en une $C(n)\epsilon^\frac{1}{384(n+1)^3}$-approximation de Hausdorff de $ M/G$ sur $\sn/A(G)$.
\end{lemma}

\begin{proof}
  Si $\sigma\in{\rm Ker}A$ alors la variété quotient de $ M$ par le groupe $<\sigma>$ vérifie aussi $\Ric\geq n-1$ et $\lambda_{n+1}\leq n+\epsilon$. Mais alors $ M$ et $ M/<\sigma>$ ont des volumes presqu'égaux à celui de $(\sn,can)$ (d'après \ref{Fiprop}). Or le cardinal de $<\sigma>$ est égal au rapport de ces volumes (car $<\sigma>$ agit librement), et donc ${\rm Ker}A=\{id_G\}$.

$G$ est un groupe agissant par isométrie sur $ M$ et $\sn$, et $\Phi$ est une application équivariante pour ces deux actions. $\Phi$ passe donc au quotient en une application de la variété riemannienne $ M/G$ sur (l'orbifold) $\sn/G$. Il est facile de montrer que $\Phi$ étant une $C(n)\epsilon^\frac{1}{384(n+1)^3}$-approxi-mation de Hausdorff de $M$ sur $\sn$ (proposition \ref{Fiprop}), son quotient est aussi une $C(n)\epsilon^\frac{1}{384(n+1)^3}$-approximation.

Pour montrer que $G$ agit librement sur $\sn$, commençons par remarquer que d'après le théorème \ref{Che-Col}, on peut supposer que sous les hypothèses du lemme, $M^n=\sn$ (remarquez que dans la suite on a juste besoin que $M$ soit homéomorphe à $\sn$, et donc \cite{Per} suffit). On utilise alors le:

\begin{lemma}\label{Suter}
  Soit $G$ un groupe fini, soient $\alpha_i:G\times\sn\to\sn$, $i=1,2$ deux actions sur la sphère $\sn$ et $\Phi:(\sn,\alpha_1)\to(\sn,\alpha_2)$ une application équivariante qui est de degré $\pm1$. Si $\alpha_1$ est libre alors $\alpha_2$ est aussi libre.
\end{lemma}

Si $\alpha_2$ n'est pas libre, alors $G$ admet un élément non trivial $\sigma$ tel que $\alpha_2(\sigma)$ fixe un point de $\sn$. Quitte à remplacer $G$ par $<\sigma^l>$ (pour $l$ bien choisi), on peut se ramener au cas où $G$ est un groupe de cardinal $p$ premier fixant un point de $\sn$ par l'action $\alpha_2$. Notons $M_\Phi$ le cylindre de l'application $\Phi$. Au couple $(M_\Phi,\sn)$ est alors associé la suite longue exacte en cohomologie relative (à coefficients dans $\mathbb{Z}/p\mathbb{Z}$) suivante:
$$\displaylines{0\to H^0(M_\Phi,\mathbb{Z}/p\mathbb{Z})\stackrel{\Phi^*}{\to}H^0(\sn,\mathbb{Z}/p\mathbb{Z})\to H^1(M_\Phi,\sn;\mathbb{Z}/p\mathbb{Z})\to\ldots\hfill\cr
\hfill\ldots\to H^n(M_\Phi,\sn;\mathbb{Z}/p\mathbb{Z})\to H^n(M_\Phi,\mathbb{Z}/p\mathbb{Z})\stackrel{\Phi^*}{\to}H^n(\sn,\mathbb{Z}/p\mathbb{Z})\to0}$$
Or, pour tout $1\leq i<n$, on a $H^i(M_\Phi,\mathbb{Z}/p\mathbb{Z})=H^i(\sn,\mathbb{Z}/p\mathbb{Z})=\{0\}$ et donc $H^i(M_\Phi,\sn;\mathbb{Z}/p\mathbb{Z})=\{0\}$. Enfin, puisque $\Phi$ est de degré $\pm1$, on en déduit que $\Phi^*$ est un isomorphisme entre $H^n(M_\Phi,\mathbb{Z}/p\mathbb{Z})$ et $H^n(\sn,\mathbb{Z}/p\mathbb{Z})$ et donc:
$$H^i(M_\Phi,\sn;\mathbb{Z}/p\mathbb{Z})=\{0\}$$
pour tout entier $i\geq0$.

Comme $\Phi$ est équivariante, les actions $\alpha_1$ et $\alpha_2$ induisent une action $\alpha_\Phi$ sur $M_\Phi$ dont les points fixes $Fix(\alpha_\Phi)$ s'identifient aux points fixes de $\alpha_2$ (car $\alpha_1$ est sans point fixe). D'après un théorème dû à E.~Floyd (théorème 7.9, chap.III, page 144 de \cite{Bre}), on a:
$$\sum_{k\geq0}{\rm dim}H^k(\underbrace{Fix(\alpha_\Phi)}_{Fix(\alpha_2)},\underbrace{Fix(\alpha_1)}_\emptyset;\mathbb{Z}/p\mathbb{Z}){\leq}\sum_{k\geq0}{\rm dim}H^k(M_\Phi,\sn;\mathbb{Z}/p\mathbb{Z}){=}0$$
et donc $H^k(Fix(\alpha_2),\emptyset;\mathbb{Z}/p\mathbb{Z})=H^k(Fix(\alpha_2);\mathbb{Z}/p\mathbb{Z})=\{0\}$ pour tout entier $k$, ce qui implique $Fix(\alpha_2)=\emptyset$. D'où une contradiction (remarquez que dans le cas qui nous interesse dans la suite, $\alpha_2$ agit par isométrie, est donc $Fix(\alpha_2)$ ne pouvait être a priori que vide ou une sous-sphère).
\end{proof}

\section{Démonstration des théorèmes B et C}
\label{sec:theoBetC}

Sous les hypothèses des théorèmes B et C, le revêtement riemannien universel (resp. orientable) $(\widetilde{M},\tilde{g})$ de la variété $ M$ vérifie $\Ric\geq n-1$ et est de volume presque égal à $\Vol\sn$. On en déduit que $(\widetilde{M},\tilde{g})$ vérifie $\lambda_{n+1}(\widetilde{M},\tilde{g})\leq n+C(n)\Bigl(\frac{\Vol\sn}{2}-\Vol M\Bigr)^{\beta(n)}$ (ce résultat découle qualitativement de \cite{Pe}; voir \cite{Au2} pour une autre preuve et le calcul de la forme du majorant). Le lemme \ref{quotient} implique alors que $\pi_1(M)$ (resp. le groupe du revêtement orientable) agit librement sur $\sn$ par isométrie et que le quotient de $\widetilde{\Phi}$ réalise une approximation de Hausdorff de $ M$ (car le revêtement est normal) sur  $\sn/\pi_1(M)$ (resp.  $\sn/(\mathbb{Z}/2\mathbb{Z})$). Pour la démonstration du théorème B, remarquons que le seul groupe d'isométries de cardinal 2 agissant librement sur $\sn$ est $\{\pm Id\}$, et donc nos variétés sont proches de $\mathbb{P}^n\Re$ en distance de Gromov-Hausdorff. Pour le théorème $C$, remarquons que le cardinal de $\pi_1(M)$ est exactement $k$ et que les groupes d'isométries de $\sn$ de cardinal $k$ sont en nombre finis (quand ils existent). On en déduit que sous les hypothèses des théorèmes B et C, les variétés sont proches d'un nombre finis de variétés riemanniennes compactes possibles de même dimension. D'après le théorème \ref{Che-Col}, on en déduit l'existence d'une constante $\alpha(n,k)>0$ telle que nos variétés soient difféomorphes aux espaces lenticulaires annoncés.

\section{Démonstration du théorème A}
\label{sec:theoA}

Nous allons en fait montrer que sur une variété vérifiant $\Ric\geq n-1$, l'existence de $n$ valeurs propres proches de $n$ implique l'existence d'une $n{+}1$-ième valeur propre proche de $n$.

Nous commençons par supposer que $M$ est orientable. On munit $E$ (comme dans la section \ref{Ln+1}) de l'orientation compatible avec celle de $TM$. Si $(S_1,\cdots,S_n)$ est une famille de sections de $E$, on note $S=S_1\wedge\cdots\wedge S_n$ la section duale de la 1-forme $S\mapsto\det(S_1,\cdots,S_n,S)$. $S_1\wedge\cdots\wedge S_n$ est bien évidemment $L^2$-orthogonale aux sections $S_i$ pour tout $i$. On a de plus les relations $D^E_X(S_1\wedge\cdots\wedge S_n)=\sum_{i=1}^nS_1\wedge\cdots\wedge D^E_XS_i\wedge\cdots\wedge S_n$ pour tout $X\in TM$ et $|S_1\wedge\ldots\wedge S_n|_E\leq|S_1|_E\ldots|S_n|_E$. On peut alors démontrer un premier résultat~:
\smallskip

\begin{lemma}\label{LnlapbrimpLn+1lapbr}
  Il existe des fonctions $\alpha(n)>0$ et $C(n)>0$ (explicitement calculables) telles que si $ M$ est orientable et vérifie $\lambda_n(\overline{\lap}^E)\leq\epsilon$ (avec $\epsilon\leq\alpha(n)$) alors $\lambda_{n+1}(\overline{\lap}^E)\leq C(n)\epsilon.$
\end{lemma}
\bigskip

\begin{proof}
  Soit $(S_1,\cdots,S_n)$ une famille $L^2$-orthonormée de sections propres de $\overline{\lap}^E$ associées aux $n$ premières valeurs propres et $S_{n+1}=S_1\wedge\cdots\wedge S_n$. Le lemme \ref{estimSi} $(iii)$ nous donne:
$$\|D^ES_{n+1}\|_2^2\leq n^2(\max_{i\leq n}\|S_i\|_\infty)^{2n-2}\max_{i\leq n}\|D^ES_i\|_2^2\leq C(n)\epsilon.$$ 
Montrons que la norme $L^2$ de $S_{n+1}$ est proche de $1$: d'après la remarque qui suit le lemme \ref{estimSi}, on a $\|S_{n+1}\|_\infty\leq\prod_{i\leq n}\|S_i\|_\infty\leq1+C(n)\sqrt{\epsilon},$
et il existe un sous-ensemble $M_\epsilon$ de $M$ tel que~:
$$\displaylines{\hfill |<S_i,S_j>_E-\delta_{ij}|\leq C(n)\epsilon^\frac{1}{4}\hfill\Vol M_\epsilon\geq\bigl(1-C(n)\epsilon^\frac{1}{4}\bigr)\Vol M\hfill}$$
Fixons un point $x$ de $M$ et $(X_1,\ldots,X_{n+1})$ un repère orthonormé direct de $E_x$ tel que ${\rm Vect}\bigl(X_1,\ldots,X_n\bigr)$ et ${\rm Vect}\bigl(S_1(x),\ldots,S_n(x)\bigr)$ coïncident, on a alors $S_1\wedge\ldots\wedge S_n=cX_{n+1}$, où $c^2$ vaut:
$$\det\begin{pmatrix}
  <S_1,S_1>_E&\ldots&<S_n,S_1>_E\\
\Bigl|&&\Bigr|\\
<S_1,S_n>_E&\ldots&<S_n,S_n>_E
\end{pmatrix}$$
On en déduit que:
$$\bigl||S_{n+1}(x)|^2-1\bigr|=|c^2-1|=\bigl|\det\bigl(<S_i,S_j>_E^{~}\bigr)_{ij}-\det I_n\bigr|\leq C(n)\epsilon^\frac{1}{4},$$
en tout point $x$ de $M_\epsilon$. On obtient donc~:
\begin{eqnarray*}
  \bigl|\|S_{n+1}\|_2^2-1\bigr|&\leq&\bigl|\frac{1}{\Vol M}\int_{M_\epsilon}\bigl(|S_{n+1}|^2-1\bigr)\bigr|\\&&+\bigl|\frac{1}{\Vol M}\int_{M\setminus M_\epsilon}\bigr(|S_{n+1}|^2-1\bigr)\bigr|\\
&\leq&C(n)\epsilon^\frac{1}{4}+\Bigl(1-\frac{\Vol M_\epsilon}{\Vol M}\Bigr)2C(n)\leq C(n)\epsilon^\frac{1}{4}
\end{eqnarray*}
$(S_i)_{1\leq i\leq n+1}$ est donc $L^2$ presque othonormée et que $\|D^ES_i\|_2^2\leq C(n)\epsilon$ pour tout $1\leq i\leq n+1$. Par l'inégalité de Cauchy-Schwarz, on obtient que pour toute combinaison linéaire $S$ des sections $S_1,\ldots,S_n$ et $S_{n+1}$ on a $\|D^ES\|_2^2\leq C(n)\epsilon\|S\|_2^2$. Le principe du min-max conclut.
\end{proof}

On peut alors en déduire la proposition suivante (remarquer qu'on ne suppose plus que $M$ est orientable) qui combinée à la proposition \ref{Fiprop} et au theorème \ref{Che-Col} démontre le théorème A~:
\smallskip

\begin{proposition}\label{LnimpLn+1}
 Il existe des constantes $\alpha(n)>0$ et $C(n)>0$ (universellement calculables) telles que si $M$ vérifie $\lambda_n(M)\leq n+\epsilon$ (avec $\epsilon\leq\alpha(n)$), alors $\lambda_{n+1}(M)\leq n+C(n)\epsilon^\frac{1}{2}.$
\end{proposition}
\medskip

\begin{remark}
  Bien évidemment, on peut décliner cette proposition en remplaçant la conclusion par l'une des inégalités~:
$$\displaylines{\Rad M\geq\bigl(1-C(n)\epsilon^{\beta(n)}\bigr)\pi,\hfill\Vol M\geq\bigl(1-C(n)\epsilon^{\beta(n)}\bigr)\Vol\S^n\cr
\mbox{ou  }d_{GH}\bigl( M,\S^n\bigr)\leq C(n)\epsilon^{\beta(n)}}$$
où $C(n)$ et $\beta(n)$ sont des constantes explicitement calculables.
\end{remark}

\begin{proof}
 Si $M$ est orientable, la proposition découle directement des propositions \ref{klapimpklapbr} $(i)$ et $(ii)$ et du lemme \ref{LnlapbrimpLn+1lapbr}. Supposons donc que $M$ n'est pas orientable. Soit $\bigl(\sqrt{n+1}.f_i\bigr)_{1\leq i\leq n}$ une famille $L^2$-orthonormée de fonctions propres associées aux $n$ premières valeurs propres non nulles de $ M$, soit $\pi:(\widetilde{M}^n,\tilde{g})\to M$ le revêtement riemannien orientable de $M$ et notons $\tilde{f}_i=f_i\circ\pi$. Sur $(\widetilde{M}^n,\tilde{g})$ on obtient une fonction $\tilde{f}_{n+1}$ telle que $(\sqrt{n+1}.\tilde{f}_i)_{i\leq n+1}$ soit une famille de fonctions propres de $(\tilde{M},\tilde{g})$ associées à des valeurs propres proche de $n$. Le morphisme de groupe $A$ de ${\rm Gr}\bigl(\widetilde{M},M\bigr)$ dans $O_{n+1}(\Re)$ construit à partir de la famille $(\tilde{f}_i)_{i\leq n+1}$ (voir le lemme \ref{quotient}) est injectif et $A\bigl({\rm Gr}\bigl(\widetilde{M},M\bigr)\bigr)$ agit librement sur $\Re^{n+1}$. Or, l'hyperplan engendré par $\bigl(\tilde{f}_i\bigr)_{1\leq i\leq n}$ est évidemment fixé par ${\rm Gr}\bigl(\widetilde{M},M\bigr)$, ce qui est contradictoire. 
\end{proof}

\begin{remark} La preuve de la proposition \ref{LnimpLn+1} peut se faire sans théorie de Flyod: on sait que l'élément non trivial de ${\rm Gr}\bigl(\widetilde{M},M\bigr)$ est d'ordre 2, fixe l'hyperplan engendré par les fonctions $(f_i\circ\pi)$ et donc agit par symétrie orthogonale hyperplane sur $\sn$. On en déduit que l'application $\widetilde{\Phi}$ de degré $\pm1$ construite à partir des fonctions propres $(\tilde{f}_i)_{i\leq n+1}$ passe au quotient en une application $\Phi$ de $M$ sur la demi-sphère $\frac{1}{2}\S^n$ dont le degré modulo $2$ est égal à 1. Or $M$ est une variété sans bord et $\frac{1}{2}\S^n$ est une variété à bord (et même contractile), le degré modulo 2 de $\Phi$ ne peut donc être que nul, d'où une contradiction.
\end{remark}

\section{Optimalité du théorème A}\label{Opt}

Pour clore cet article, nous allons montrer que pour tout entier $p\in\{1,\ldots,n{-}1\}$, il existe une suite $(g_k)_{k\in\N}$ de métriques sur $\S^n$ qui vérifie $\Ric\geq(n{-}1)$, $\lambda_{p}(g_k)\to n$, $\lambda_{p+1}(g_k)\to \beta(n)>n$, $Vol(g_k)\to 0$, $(\S^n,g_k)$ tend en distance de Gromov-Hausdorff vers la demi-sphère canonique de dimension $p$ et par conséquent, la suite des Radius tend vers $\frac{\pi}{2}$ (ce qui montre qu'aucun des résultats de stabilité obtenu dans la proposition \ref{LnimpLn+1} et dans la remarque qui la suit ne se généralise sous l'hypothèse $\lambda_{n-1}\leq n+\epsilon$). En revanche, nous n'avons pas encore réussi à construire une suite de variétés non difféomorphes à $\sn$ de courbure de Ricci supérieure à $(n{-}1)$ et telle que la suite des $\lambda_{n-1}$ tende vers $n$ (voir \cite{An2} ou \cite{Ots} pour des variétés non homotopes à $\S^n$ admettant des métriques de courbure de Ricci supérieure à $(n{-}1)$ et admettant une seule valeur propre arbitrairement proche de $n$).
\smallskip

\noindent{\sc{\bf Exemple.}} --- La métrique $g_k$ est une généralisation des fuseaux en codimension quelconque. En fait, on construit $g_k$ en écrasant la métri-que canonique de $\sn$ dans le fibré normal à une sous sphère $\S^{p-1}$. La difficulté étant de régulariser la métrique sur $\S^{p-1}$ et sur son cut-locus $\S^{n-p}$ tout en préservant le minorant sur la courbure de Ricci. Nous décrivons la suite $g_k$ dans le cas particulier $p=n-1$, la généralisation aux autres valeurs de $p$ ne pose aucun problème (notons que dans le cas $p=n$, il n'est pas possible d'écraser la métrique dans le facteur normal).
 Soit $k$ un entier non nul. Nous commençons par fixer les notations suivantes:
$$\displaylines{\hfill\eta_k=\sqrt{\sin^2(\frac{1}{\sqrt{k}})+\frac{1}{k^2}\cos^2(\frac{1}{\sqrt{k}})},\hfill\epsilon_k=\frac{\frac{\pi}{2}-\frac{1}{\sqrt{k}}}{k},\hfill\cr
\theta_k=\arctan\Bigl(\frac{1}{k\tan\frac{1}{\sqrt{k}}}\Bigr)-\frac{\pi}{2k}+\frac{1}{k\sqrt{k}}}$$
On note $I_k$ lintervalle $]0,\frac{\pi}{2}-\theta_k[$. Sur la variété $M=I_k\times\S^1\times\S^{n-2}$, on considère la métrique $g_k=dr^2+a_k(r)^2g_{\S^1}+b_k(r)^2g_{\S^{n-2}}$, où $g_{\S^1}$ et $g_{\S^{n-2}}$ sont les métriques canoniques des sphères $\S^1$ et $\S^{n-2}$, et où $a_k$ et $b_k$ sont des fonctions définies sur $I_k$ par les formules~:
$$\displaylines{a_k(r)=\left\{
  \begin{array}{ll}
\frac{1}{k}\sin(kr)&\mbox{sur }]0,\epsilon_k],\\
\eta_k\sin(r+\theta_k)&\mbox{sur }[\epsilon_k,\frac{\pi}{2}-\theta_k[,\end{array}\right.\hfill\cr
\mbox{ et }\cr
 b_k(r)=\left\{
  \begin{array}{ll}
\frac{\epsilon_k}{\epsilon_k+\theta_k}\cos\bigl((\frac{\epsilon_k+\theta_k}{\epsilon_k})r\bigr)+\frac{\theta_k}{\epsilon_k+\theta_k}\cos\bigl(\theta_k+\epsilon_k\bigr)&\mbox{sur }]0,\epsilon_k],\\
\cos(r+\theta_k)&\mbox{sur }[\epsilon_k,\frac{\pi}{2}-\theta_k[.
  \end{array}\right.\hfill}$$
Remarquons que la fonction $a_k$  (resp. $b_k$) est $C^1$ sur $I_k$, $C^2$ en dehors de $\epsilon_k$ et tend vers $0$ en $0$ (resp. tend vers $0$ en $\frac{\pi}{2}-\theta_k$). La métrique $g_k$ se prolonge en une métrique $C^1$ sur $\S^n$ pour laquelle le cercle $\S^1$ est le cut-locus de la sous sphère $\S^{n-2}$.

Pour minorer la courbure de Ricci, on applique les formules de calcul de la courbure de Ricci des doubles-produits tordus (voir par exemple \cite{Ots}). En un point $x$ fixé de $M$ on note $\frac{\partial}{\partial r}$ un vecteur tangent au facteur $I_k$, $u$ un élément de $T_x\S^1$, $v$ un élément de $T_x\S^{n-2}$. Alors, la courbure de Ricci $\Ric_k$ de $(\sn,g_k)$ vérifie les propriétés suivantes:

$$\Ric_k\bigl(\frac{\partial}{\partial r},u\bigr)=\Ric_k\bigl(\frac{\partial}{\partial r},v\bigr)=\Ric_k\bigl(u,v\bigr)=0,$$
$$\displaylines{\Ric_k\bigl(\frac{\partial}{\partial r},\frac{\partial}{\partial r}\bigr)=-\frac{a''}{a}-(n-2)\frac{b''}{b}\hfill\cr
\hfill=\left\{
  \begin{array}{l}
k^2+(n-2)\bigl(\frac{\epsilon_k+\theta_k}{\epsilon_k}\bigr)^2\frac{\cos\bigl((\frac{\epsilon_k+\theta_k}{\epsilon_k})r\bigr)}{\cos\bigl((\frac{\epsilon_k+\theta_k}{\epsilon_k})r\bigr)+\frac{\theta_k}{\epsilon_k}\cos\bigl(\theta_k+\epsilon_k\bigr)}\geq k^2\mbox{ sur }]0,\epsilon_k]\\
n-1\mbox{ sur }[\epsilon_k,\frac{\pi}{2}-\theta_k[\\
  \end{array}\right.}$$
$$\displaylines{\frac{\Ric_k\bigl(u,u\bigr)}{g_k(u,u)}=-\frac{a''}{a}-(n-2)\frac{a'b'}{ab}\hfill\cr
\hfill=\left\{
  \begin{array}{l}
k^2+(n-2)\frac{k\cos(kr)}{\sin(kr)}\bigl(\frac{\epsilon_k+\theta_k}{\epsilon_k}\bigr)\frac{\sin\bigl((\frac{\epsilon_k+\theta_k}{\epsilon_k})r\bigr)}{\cos\bigl((\frac{\epsilon_k+\theta_k}{\epsilon_k})r\bigr)+\frac{\theta_k}{\epsilon_k}\cos\bigl(\theta_k+\epsilon_k\bigr)}\geq k^2\mbox{ sur }]0,\epsilon_k]\\
n-1\mbox{ sur }[\epsilon_k,\frac{\pi}{2}-\theta_k[\\
  \end{array}\right.}$$
et 
$$\frac{\Ric_k\bigl(v,v\bigr)}{g_k(v,v)}=-\frac{b''}{b}-\frac{a'b'}{ab}+(n-3)\Bigl(\frac{1-{b'}^2}{b^2}\Bigr)$$
Donc
$$\displaylines{\frac{\Ric_k\bigl(v,v\bigr)}{g_k(v,v)}\geq\Bigl(\frac{\epsilon_k+\theta_k}{\epsilon_k}\Bigr)^2\frac{\cos\bigl((\frac{\epsilon_k+\theta_k}{\epsilon_k})r\bigr)}{\cos\bigl((\frac{\epsilon_k+\theta_k}{\epsilon_k})r\bigr)+\frac{\theta_k}{\epsilon_k}\cos\bigl(\theta_k+\epsilon_k\bigr)}\hfill\cr
\hskip1.7cm\geq\frac{\epsilon_k+\theta_k}{\epsilon_k}\sim\frac{2}{\pi}\sqrt{k}\mbox{ sur }]0,\epsilon_k]\hfill}$$
(car $\theta_k+\epsilon_k\geq\bigl(\frac{\epsilon_k+\theta_k}{\epsilon_k}\bigr)r$) et
$$\displaylines{\frac{\Ric_k\bigl(v,v\bigr)}{g_k(v,v)}=
2+(n-3)\frac{1-\sin^2(r+\theta_k)}{\cos^2(r+\theta_k)}=n-1\mbox{ sur }[\epsilon_k,\frac{\pi}{2}-\theta_k[.}$$
On déduit des calculs précédents que $(\S^n,g_k)$ est de courbure de Ricci supérieure à $(n{-}1)$ pour $k$ assez grand.

\noindent Il est évident que le volumes des métriques $g_k$ tend vers $0$. De plus, la suite de variétés ainsi obtenue tend (en distance de Gromov-Hausdorff) vers une des hémi-sphères de dimension $n{-}1$ que borde la sphère $\S^{n-2}$: notons $(N_k,h_k)$ la variété $]\epsilon_k,\frac{\pi}{2}-\theta_k[\times\S^{n-2}$, munie de la métrique $h_k=(dr)^2+\cos^2(r+\theta_k)\,g_{\S^{n-2}}$; elle est isométrique à une boule géodésique de $\S^{n-1}$ de rayon $\frac{\pi}{2}-(\epsilon_k+\theta_k)$ (privée de son centre), donc elle converge vers la demi-sphère $\frac{1}{2}\S^{n-1}$, munie de sa métrique canonique. Par ailleurs, $(N_k,h_k)$ se plonge dans $(M,g_k)$ de manière isométrique, via l'application $(r,v)\mapsto(r,u_0,v)$, où $u_0$ est un point fixé de $\S^1$. De plus, la distance dans $M$ entre deux points $p=(r,u_0,v)$ et $q=(r',u_0,v')$ coïncide avec la distance dans $N_k$ entre $(r,v)$ et $(r',v')$. Comme $d_{g_k}\bigl[(r,u,v);(r,u_0,v)\bigr]\leq\pi\,\eta_k$, on obtient que $(M,g_k)$ converge vers l'hémisphère de $\S^{n-1}$ au sens de Gromov-Hausdorff. Il nous reste à montrer que le laplacien de ces variétés admet au moins $n{-}1$ valeurs propres proches de $n$. Pour cela, nous pouvons soit appliquer \cite{Ber}, soit choisir des points $x_0,\ldots,x_{n-2}$ sur $\S^{n-2}$ de sorte que $d_{\S^{n-2}}(x_i,x_j)=\frac{\pi}{2}$ si $i\neq j$ et noter encore $x_i$ le point $(0,u_0,x_i)$ de $M$; en appliquant le principe du min-max aux $(n-1)$ fonctions $f_i(x)=\cos\bigl(d_k(x_i,x)\bigr)$ (la forme des métriques $g_k$ fait que les intégrales à calculer convergent vers les intégrales correspondantes calculées sur $(\frac{1}{2}\S^{n-1},can)$; voir \cite{Au2} pour les détails). Enfin, les métriques $g_k$ n'ont  pas plus de $n-1$ valeurs propres proches de $n$ car sinon, d'après \cite{Ber}, les variétés $(\sn,g_k)$ contiendraient une sous partie Gromov-Hausdorff proche de $(\S^{n-1},can)$; or une demi-sphère de dimension $n-1$ ne peut contenir une partie Gromov-Hausdorff proche de la sphère $(\S^{n-1},can).$ On remarquera toutefois que la $p$-ième valeur propre $\lambda_p^k$ de la variété $(M,g_k)$ reste bornée par une fonction de $p$ (pour $p$ quelconque) lorsque $k$ tend vers l'infini car les fonctions propres radiales pour le problème de Neuman de la demi-sphère de dimension $n{-}1$ permettent de construire, comme précédemment, une famille presque $L^2$-orthonormée de fonctions tests sur $M$ dont les quotients de Rayleigh pour la métrique $g_k$ tendent vers le spectre pour le problème de Neuman de la demi-sphère de dimension $n{-}1$. On conclut alors par le principe du min-max.


\begin{thebibliography}{aa}
\bibitem{An2} {\sc M.~Anderson}, {\em Metrics of positive Ricci curvature with large diameter}, Manuscripta Math. {\bf 68} (1990), p. 405--415.

\bibitem{Au1} {\sc E.~Aubry}, {\em Théorème de la Sphère}, Séminaire de théorie spectrale et géométrie, Grenoble Volume {\bf 18} (2000), p. 125--155.

\bibitem{Au2} {\sc E.~Aubry}, {\em Variétés de courbure de Ricci presque minorée : inégalités géométriques optimales et stabilité des variétés extrémales}, Thèse, Institut Fourier, Grenoble (2003).

\bibitem{BBG} {\sc P.~Bérard, G.~Besson, S.~Gallot}, {\em Sur une inégalité isopérimétrique qui généralise celle de Paul Lévy-Gromov},  Invent. Math. {\bf 80} (1985), p. 295--308.

\bibitem{Ber} {\sc J.~Bertrand}, Preprint, Université Paris-Sud, Orsay  (2004).

\bibitem{Bre} {\sc G.~Bredon}, {\em Introduction to compact transformation groups}, Pure and Applied Mathematics, Vol {\bf 46} Academic Press New-york--London (1972).

\bibitem{Ch-Co3} {\sc J.~Cheeger,~T.~Colding}, {\em On the structure of spaces with Ricci curvature bounded below. I} J. Differential Geom. {\bf 46} (1997), p. 406--480.

\bibitem{Co1} {\sc T.~Colding},  {\em Shape of manifolds with positive Ricci curvature}, Invent. Math. {\bf 124} (1996), p. 175--191.

\bibitem{Co2} {\sc T.~Colding},  {\em Large manifolds with positive Ricci curvature}, Invent. Math. {\bf 124} (1996), p. 193--214.

\bibitem{Ga} {\sc S.~Gallot}, {\em Volume, courbure de Ricci et convergence des variétés (d'après T. H. Colding et Cheeger-Colding)}, séminaire Bourbaki Nov 1997 (exposé n°835), Astérisque {\bf 252} (1998), p. 7--32.

\bibitem{Ili} {\sc S.~Ilias}, {\em Constantes explicites pour les inégalités de Sobolev sur les variétés riemanniennes compactes}, Ann. Inst. Fourier {\bf 33} (1983), p.151-165.

\bibitem{Ots} {\sc Y.~Otsu},  {\em On manifolds of positive Ricci curvature with large diameter}, Math. Z. {\bf 206} (1991), p. 252--264.

\bibitem{Pe} {\sc P.~Petersen},  {\em On eigenvalue pinching in positive Ricci curvature}, Invent. Math. {\bf 138} (1999), p. 1--21.

\bibitem{Pe2} {\sc P.~Petersen},  {\em On eigenvalue pinching in positive Ricci curvature, Erratum}, Invent. Math. {\bf 155} (2004), p. 223.

\bibitem{Per} {\sc G.~Perelman},  {\em Manifold of positive Ricci curvature with almost maximal volume} JAMS {\bf 7} (1994), p. 299--305.

\bibitem{Ru} {\sc E. Ruh}, {\em Curvature and differentiable structure on spheres}, Bull. Amer. Math. Soc. {\bf 77} (1971), p. 148--150.

\bibitem{Yam2} {\sc T. Yamaguchi}, {\em Lipschitz convergence of manifolds of positive Ricci curvature with large volume}, Math. Ann. {\bf 284} (1989), p. 423--436.
\end{thebibliography}
\end{document}